\documentclass[epjST]{svjour}
\usepackage{hyperref}
\usepackage{graphics}

\usepackage{dsfont}

\usepackage{amsmath}
\usepackage{amssymb}
\usepackage{subfigure}

\def\mydate{\number\day\ {\ifcase\month \or January\or February\or
              March\or April\or May\or June\or July\or August\or
              September\or October\or November\or December\fi}
\number\year}

\def \balpha{{\boldsymbol\alpha}}
\def \fb{\mathbf{f}}
\def \nb{\mathbf{n}}
\def \bthe{{\boldsymbol \Theta}}
\def \bmu{{\boldsymbol \mu}}
\def \OmegaI{\Omega_{\mathcal{I}}}

\newcommand{\mbR}{\mathbb{R}}

\newcommand{\sint}[4]{\int_{#1}^{#2}{#3}\,d #4}
\newcommand{\ds}{\displaystyle}
\newcommand{\mcD}{\mathcal{D}}

\begin{document}

\title{The exit-time problem for a Markov jump process}

\author{Nathanial Burch\inst{1}\fnmsep\thanks{\email{burchn@gonzaga.edu}} Marta D'Elia\inst{2}\fnmsep\thanks{\email{mdelia@sandia.gov}} \and R. B. Lehoucq\inst{2}\fnmsep\thanks{\email{rblehou@sandia.gov}} }
\institute{Department of Mathematics, Gonzaga University, Spokane, WA 99258 \and Sandia National Laboratories, Albuquerque, NM 87185-\{1321, 1320\}. Sandia is a multi-program laboratory managed and operated by Sandia Corporation, a wholly subsidiary of Lockheed Martin Corporation, for the U.S. Department of Energy's National Nuclear Security Administration under contract DE-AC04-94AL85000.}

\abstract{
%
The purpose of this paper is to consider the exit-time problem for a finite-range Markov jump process, i.e, the distance the particle can jump is bounded independent of its location.
Such jump diffusions are expedient models for anomalous transport exhibiting super-diffusion or nonstandard normal diffusion. 
We refer to the associated deterministic
equation as a volume-constrained nonlocal diffusion equation. The volume constraint is the nonlocal
analogue of a boundary condition necessary to demonstrate that the nonlocal diffusion equation is well-posed and is consistent with the jump process. A critical aspect of the analysis is a variational formulation and a recently developed nonlocal vector calculus. This calculus allows us to pose nonlocal backward and forward Kolmogorov equations, the former equation granting the various moments of the exit-time distribution.
} 
\maketitle
\begin{center}
\emph{Draft as of} \mydate
\end{center}

\section{Introduction}
\label{intro}

The classical Brownian motion model for diffusion is not well-suited for applications with discontinuous sample paths. For instance, the mean square displacement of a diffusing particle associated with a jump process often grows faster than that for the case of Brownian motion, or grows at the same rate but is of finite variation or activity (terms that we will define in \S\ref{sec:fluc}).  Such jump diffusions are expedient models for anomalous transport exhibiting super-diffusion or nonstandard normal diffusion. Examples include various problems in finance, transport in heterogenous media~\cite{neta:09}, species migration or heterogenous heat conduction as suggested by the length and time scales over which the data is collected. See the edited volume for a collection of papers \cite{klrs} containing further information and an abundance of references.

The purpose of this manuscript is to discuss the exit-time problem for finite-range Markov jump processes. A finite-range process restricts the jump-rate to be zero outside of a bounded neighborhood about the location of the particle. Because the process may jump outside of the domain, the associated deterministic equation, in contrast to the classical Fokker-Planck, must be augmented with a volume-constraint instead of a boundary condition. The volume-constraint is the nonlocal analogue of a boundary condition necessary to demonstrate that the equation is well-posed and is consistent with the jump process. In particular, realizing a Monte Carlo simulation to compute various exit-time statistics mandates absorbing (or killing) the particle upon departure out of the domain that rarely, if ever, occurs at the boundary. In fact, with probability one, a particle originating in a domain does not touch the boundary; see, e.g., \cite{mill:75}. Hence enforcing a volume-constraint for the deterministic equation is symbiotic with Monte Carlo simulations. 

A distinguishing aspect of our treatment for the exit-time problem is that the deterministic equation is given by a master equation. This allows us to analyze spatially inhomogeneous problems where the jump-rate depends upon location and may be asymmetric, i.e., the rate to and from a point may be distinct.
A critical aspect is a recently developed nonlocal vector calculus that enables striking analogies to be drawn with the classical vector calculus including Fick's laws and the backward, forward Kolmogorov equations requiring the notion of an adjoint operator. The flexibility afforded by the nonlocal vector calculus enables us to consider exit-time problems over nontrivial domains in $\mbR^n$ and builds upon work accomplished via the use of fractional derivative based approaches; see, e.g., \cite{delc:06,cade:07,chmn:12,kdmm:13} and the references provided. Although we have assumed that the jump-rate is that associated with a finite-range jump process, the representation of the master equation in terms of the calculus is valid for a multitude of classes of jump-rates, including those of infinite range such as L\'{e}vy measures or their tempered, truncated variants; see \cite{mast:94,cade:07,bame:10} including the recent review \cite{samc:14}. Said another way, the nonlocal vector calculus is jump-rate agnostic and oblivious to a finite domain.

The nonlocal vector calculus also lays the foundation for a variational formulation of the deterministic exit-time problem.  We can then establish that a broad range of volume-constrained problems are well-posed. This lays the foundation for stable numerical methods of the volume-constrained problem. This calculus also enables us to pose the nonlocal analogues of the forward and backward Kolmogorov equations, useful for determining the various moments of the exit-time distribution. 

The exit-time problem for L\'{e}vy motion has received attention in the research literature, see, e.g., \cite{gdls:14} for a recent reference and in particular, Brownian motion is well-studied. In comparison, the general case of a spatially inhomogenous Markov process is treated sporadically in the literature; for instance \cite{wesz:83} considers the one-dimensional problem, and the paper \cite{kmst:86} explains that a ``boundary layer'', or what amounts to a volume-constraint, is needed in addition to a boundary condition. To the best of our knowledge, our paper is the first general treatment for the deterministic problem associated with the exit-time distribution for a general class of Markov ``pure'' jump processes, i.e., no Brownian component or killing term.


\section{Markov jump process} \label{markov-jump-proess}

The master equation
\begin{align}
 u_t(x,t)  & =  \sint{\mbR^n}{}{\big( u(y,t)\gamma(y,x) - u(x,t)\gamma(x,y) \big)}{y} \qquad x \in \mbR^n\,, \label{fr-sp-ms-eq}
 \end{align}
is useful for describing a spatially inhomogenous Markov process; see, e.g., \cite{klso:11} for a discussion. The two-point function $ \gamma: \mathbb{R}^n \times \mathbb{R}^n \to \mathbb{R}$ and $u:\mathbb{R}^n \to \mathbb{R}$ are the jump-rate and the density of particles, respectively.
The associated Markov jump process can be simulated by Monte Carlo realizations of a continuous-time Markov chain over a continuum state-space, or equivalently, as realizations of an off-lattice continuous-time random walk (CTRW) with an exponential jump-rate. The realizations are constructed from the master equation \eqref{fr-sp-ms-eq}. The equation explains that the temporal rate of change of probability $u(x,t)\,dx$ of locating a particle at $x$ about the volume $dx$ at time $t$ is given by the difference in probability gain and loss at the point $x$. The probability gain is given by the jump-rate $\gamma(y,x)\,dx$ into $dx$ from $y$ given the probability $u(y,t)\,dy$ and in analogous fashion, the probability loss is given by the jump-rate $\gamma(x,y)\,dy$ into $dy$ from $x$  given the probability $u(x,t)\,dx$. As we will demonstrate in \S\ref{sec:ncde}, the master equation embodies a nonlocal Fick's first and second law of diffusion.

We also further suppose that the jump-rate $\gamma(x,y)$ is zero outside a ball of radius $\lambda$, i.e.,
\begin{align}\label{assump-jr}
\gamma(x,y) & = 0 \,, \text{ for } |x-y| \geq \lambda>0\,.
\end{align}
In words, the particle can jump to $y$ from $x$ when the distance between $x$ and $y$ is no larger than $\lambda$.
We say that ``nonlocal convection'' occurs when $\gamma(x,y) \neq \gamma(y,x)$. The fluctuations of the particle path can take a myriad of forms; see \S\ref{sec:fluc} for a discussion. The scaling of the mean square displacement for the jump diffusion ultimately depends upon the analytic properties of the jump-rate. The scaling is either linear or larger so that the diffusion is either normal or superdiffusive. Truncating a radial L\'{e}vy measure removes the ``heavy tails'' rendering a normal jump diffusion that is not a Brownian motion; see \cite{mast:94}.

The master equation \eqref{fr-sp-ms-eq} is an instance of the more general diffusion equation
\begin{align}
	u_t(x,t) & = \sint{\mbR^n}{}{\big( h(y,x,t) - h(x,y,t) \big)}{y}\,, \notag
\end{align}
via the relationship $h(x,y,t)=u(x,t)\gamma(x,y)$.
Such an equation is a nonlocal analogue of $ u_t  = -\nabla \cdot \mathbf{q}$. The nonlocal diffusion arises because points $y\neq x$ determine the rate of diffusion at $x$.
An expression for $h$ (or $\mathbf{q}$) can be derived as an ensemble average in phase space; see the paper \cite{lese:11}.
The conclusion is that there is a basis for nonlocal diffusion in non-equilibrium statistical mechanics; the classical diffusion arises by postulating Fick's first law, or assuming that the particle sample path is continuous. In particular, \S\ref{sec:ncde} postulates a nonlocal analogue of Fick's first law using the nonlocal vector calculus reviewed in \S\ref{sec:NLVC}.

\section{Finite domain} \label{sec:process-cnfd-dm}

\begin{figure}[htbp]
	\begin{center}
    \resizebox{0.75\columnwidth}{!}{\includegraphics{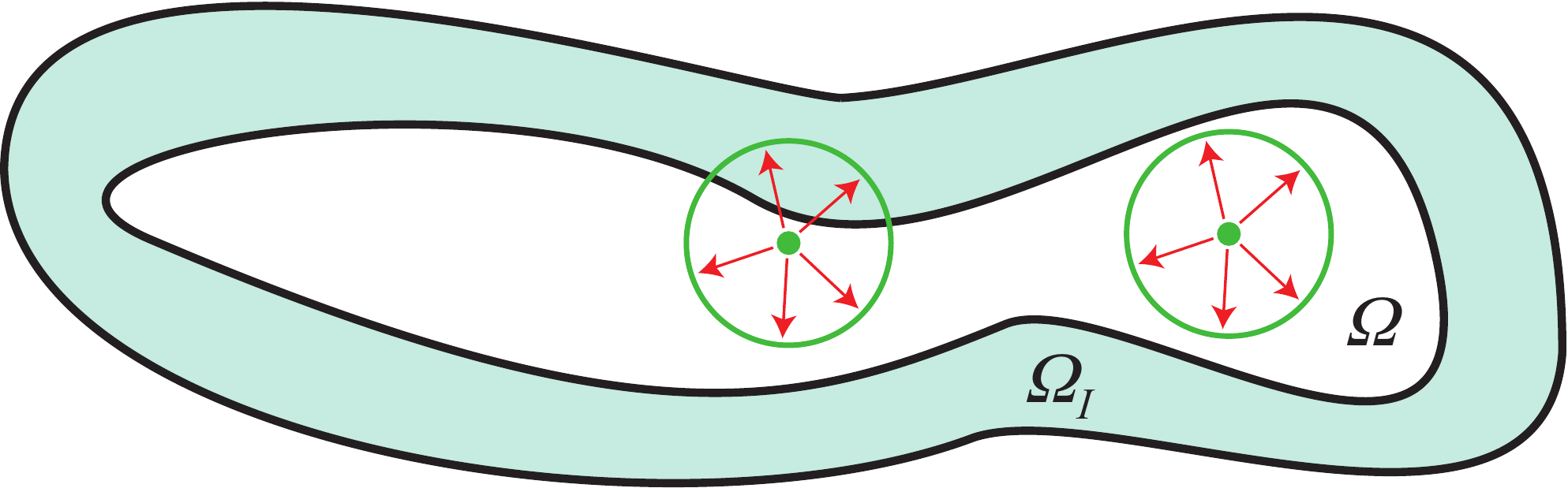}}
	\end{center}
\caption{Finite domain $\Omega $ and interaction domain $\OmegaI$ given a finite-range process; the particle can jump within the ball about the point located at center.} \label{fig:vc}
\end{figure}

\begin{figure}[htbp]
	\begin{center}
	\resizebox{0.75\columnwidth}{!}{\includegraphics{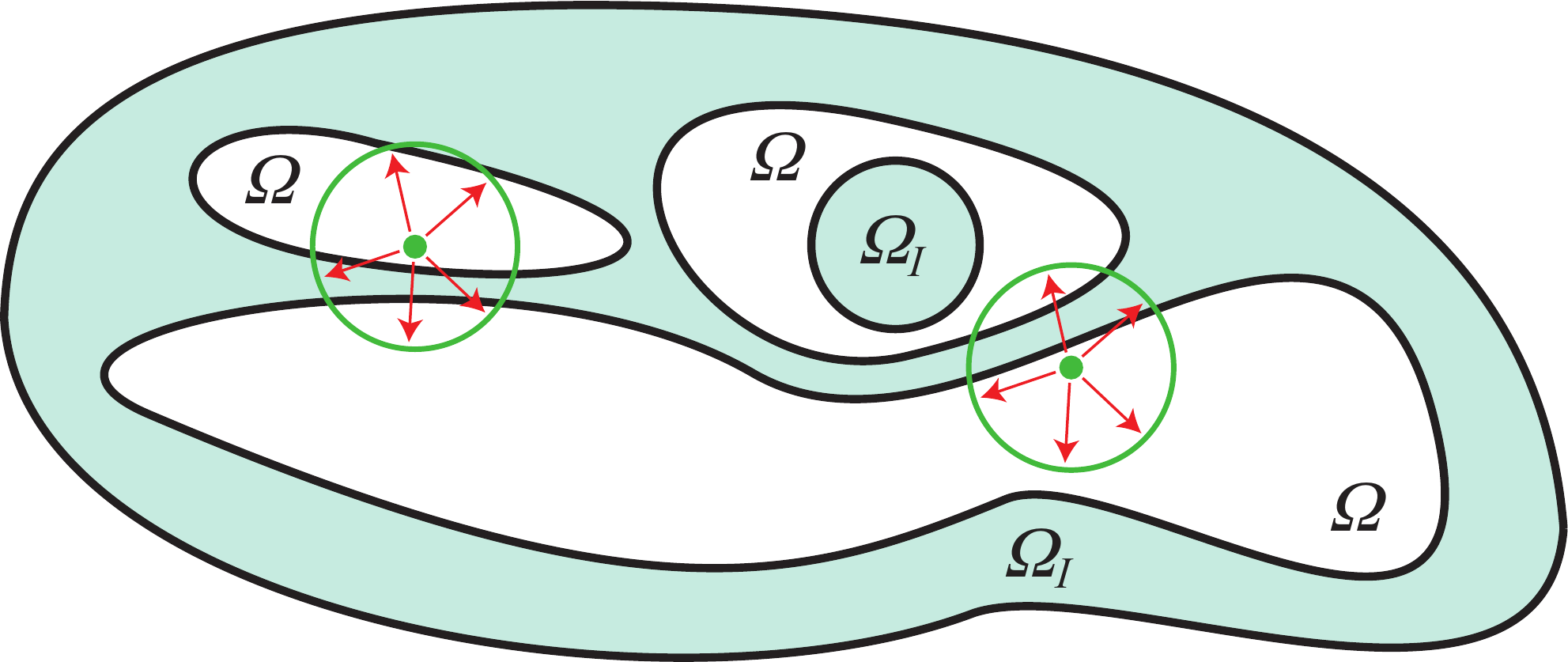}}
	\end{center}
\caption{A more elaborate example of a domain $\Omega $ and its interaction domain $\OmegaI$ given a finite-range process; the particle can jump within the ball about the point located at center.} \label{fig:vc2}
\end{figure}

Suppose we have a finite domain, $\Omega $, of interest. By assumption \eqref{assump-jr}, the jump-rate $\gamma$ is restricted so that the particle can jump at most a finite distance $\lambda$ from $x$. We define the interaction domain $\OmegaI$ to be the region to which the particle may jump to when originating in $\Omega$. Hence, by our assumption \eqref{assump-jr} on the jump-rate, the interaction domain is also finite. Figure \ref{fig:vc} displays an example of the interaction domain
\[
\OmegaI = \bigcup_{x\in\Omega} B_\lambda(x) \setminus \Omega\,,
\]
where $B_\lambda(x)$ denotes the ball about $x$ of radius $\lambda$. Such an interaction domain provides a collar for the domain $\Omega$. Figure \ref{fig:vc2} provides a more elaborate example of a disconnected domain and its interaction domain.

The master equation \eqref{fr-sp-ms-eq} describes a spatially inhomogeneous system undergoing a possibly asymmetric jump-rate (or nonlocal convection) over all of $\mbR^n$. An example of such a system is given by the jump-rate kernel
\begin{align}
   \mathds{1}_{\Omega\cup \OmegaI}(x)\, \mathds{1}_{\Omega\cup \OmegaI}(y) \nu(x-y) \label{ex-jr}
\end{align}
where $\nu$ is a L\'{e}vy measure and
$\mathds{1}_{\Omega\cup \OmegaI}$ is the indicator function given by
\begin{equation} \label{indicator}
\mathds{1}_{\Omega\cup \OmegaI}(x) :=
\begin{cases}
1 & x \in \Omega\cup\OmegaI \,,\\
0 & x \notin \Omega\cup\OmegaI\,.
\end{cases}
\end{equation}
Inserting the jump-rate kernel \eqref{ex-jr} into the free-space master equation \eqref{fr-sp-ms-eq} leads to the equation
\begin{subequations}
\begin{align}
	u_t(x,t) & =  \sint{\Omega \cup \OmegaI}{}{\big( u(y,t)\nu(y) - u(x,t)\nu(x) \big)}{y} \qquad x \in \mbR^n\,. \notag
\end{align}
We remark that the assumption \eqref{assump-jr} on the jump-rate kernel satisfies
\begin{align}
\gamma(x,y)& =0  \text{ for } x\in \Omega  \text{ and } y \in \mathbb{R}^n\setminus\big(\Omega\cup \OmegaI\big)\,, \label{assump-gamma}
\end{align}
and leads to the confined nonlocal diffusion equation
\begin{align}
	u_t(x,t) & =  \sint{\Omega \cup \OmegaI}{}{\big( u(y,t)\gamma(y,x) - u(x,t)\gamma(x,y) \big)}{y} \qquad x \in \Omega\cup\OmegaI\,. \label{conf-ms-eq}
\end{align}
\end{subequations}
Realizations of the Markov jump process corresponding to the confined master equation \eqref{conf-ms-eq} are constructed as described in the discussion following \eqref{fr-sp-ms-eq}.

\begin{figure}[htbp]
	\begin{center}
    \resizebox{0.75\columnwidth}{!}{\includegraphics{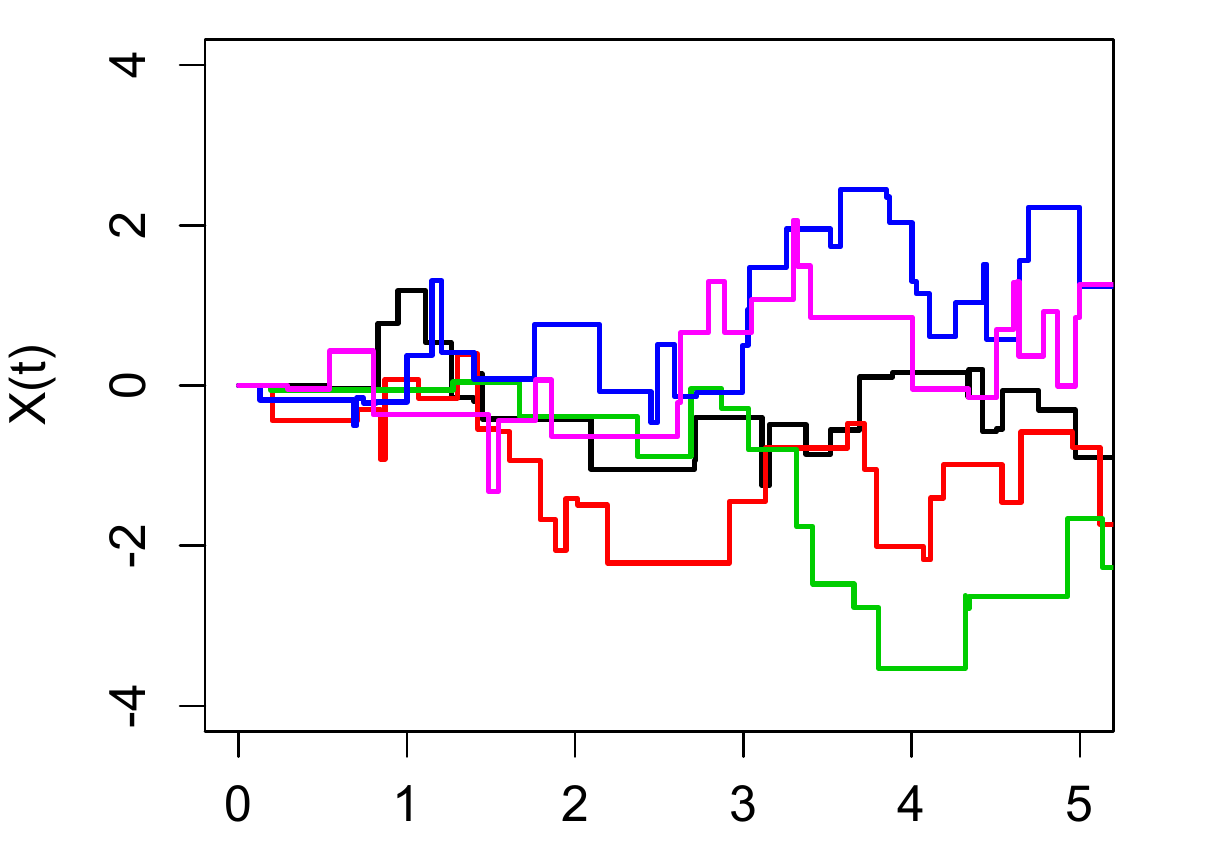}}
	\end{center}
\caption{Realizations of a compound Poisson process, a finite-activity, finite variation process; wait-times are exponentially distributed with parameter set to $0.2$. The jumps are uniformly distributed in the interval $(-1, 1)$ so that the process is of finite range.} \label{fig:sp1}
\end{figure}
\begin{figure}[htbp]
	\begin{center}
    \resizebox{0.75\columnwidth}{!}{\includegraphics{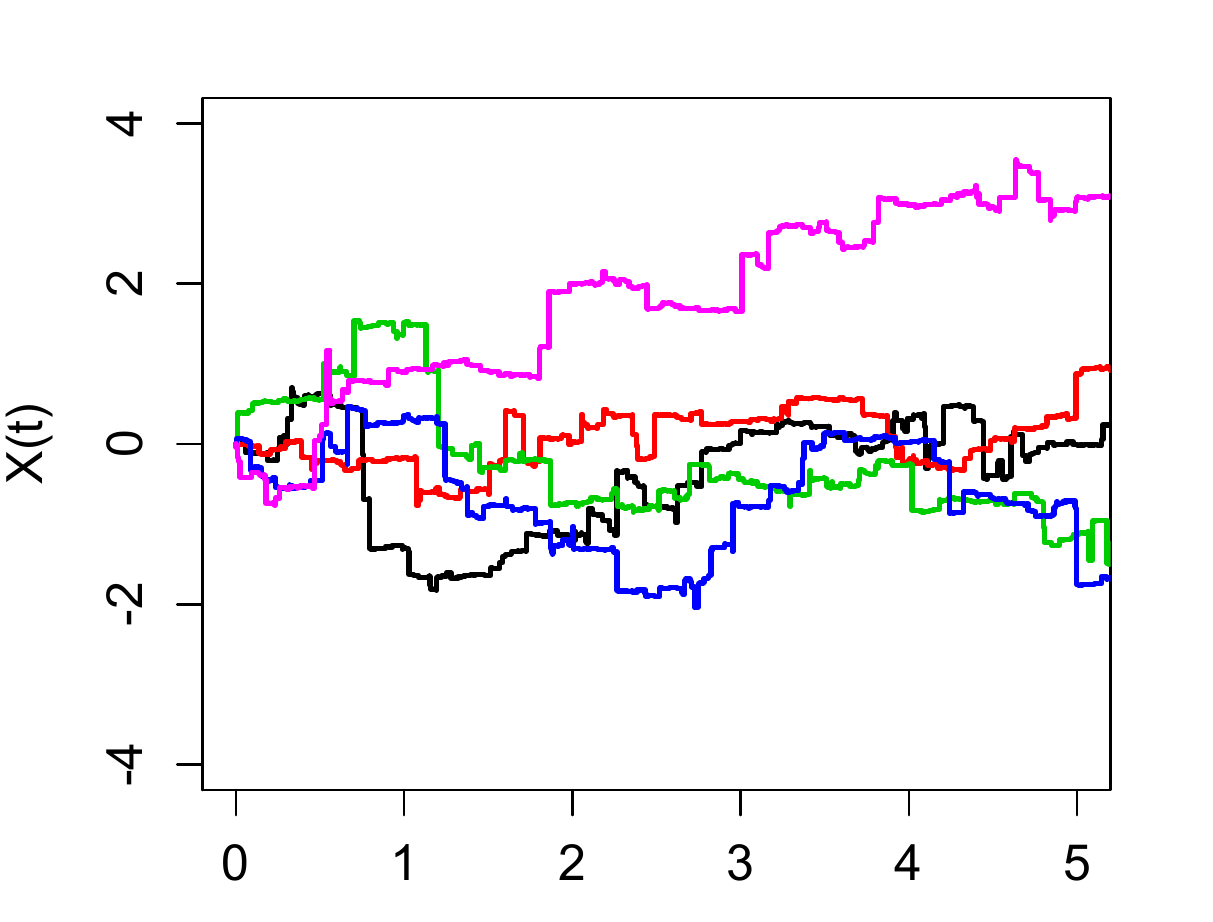}}
	\end{center}
\caption{Realizations of a infinite activity, finite variation process with jump-rate kernel \eqref{nu} approximated by \eqref{nu-app} with $\alpha = 1/2$, $\lambda = 1$, $m = 1$, $\varepsilon = 0.001$.} \label{fig:sp2}
\end{figure}
\begin{figure}[htbp]
	\begin{center}
    \resizebox{0.75\columnwidth}{!}{\includegraphics{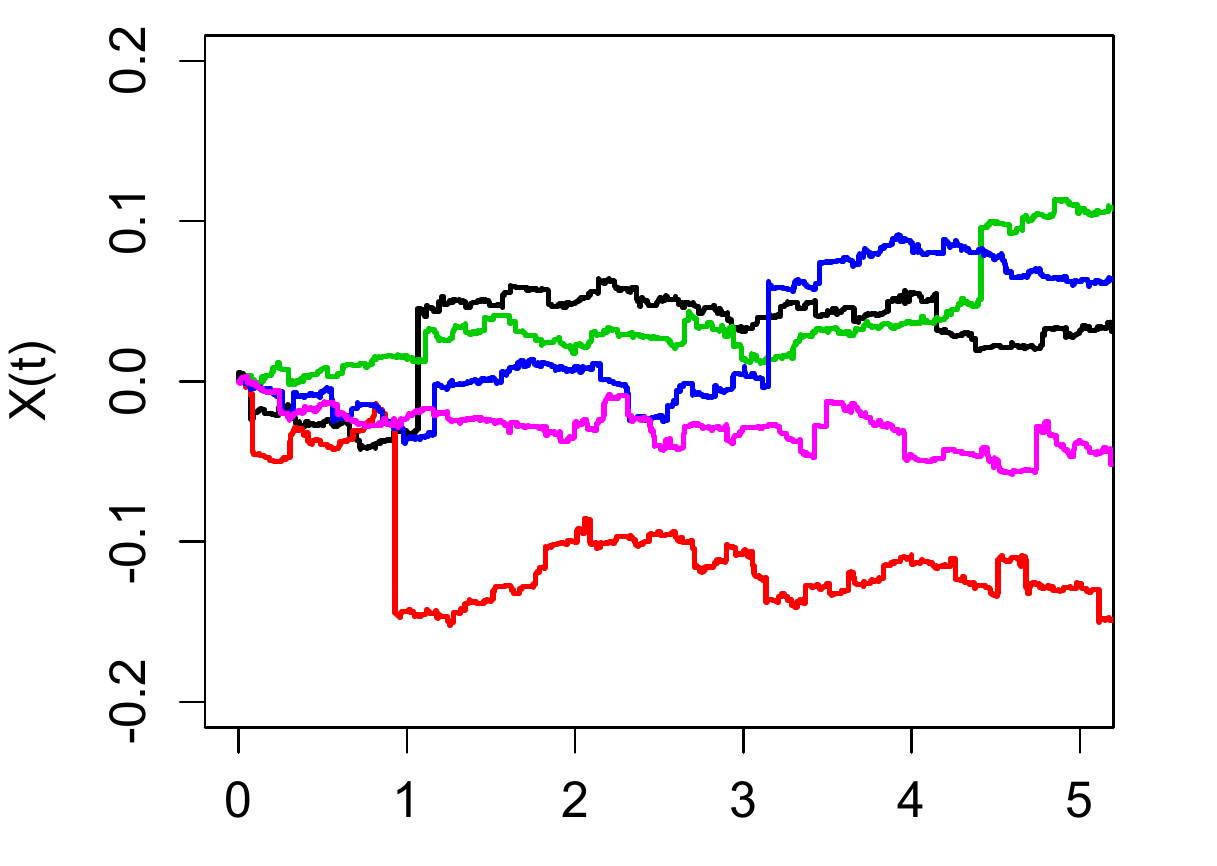}}
	\end{center}
\caption{Realizations of a infinite activity, infinite variation process  with jump-rate kernel \eqref{nu} approximated by \eqref{nu-app} with $\alpha = 3/2$, $\lambda = 1$, $m = 1000$, and $\varepsilon = 0.001$.} \label{fig:sp3}
\end{figure}
\begin{figure}[htbp]
	\begin{center}
    \resizebox{0.75\columnwidth}{!}{\includegraphics{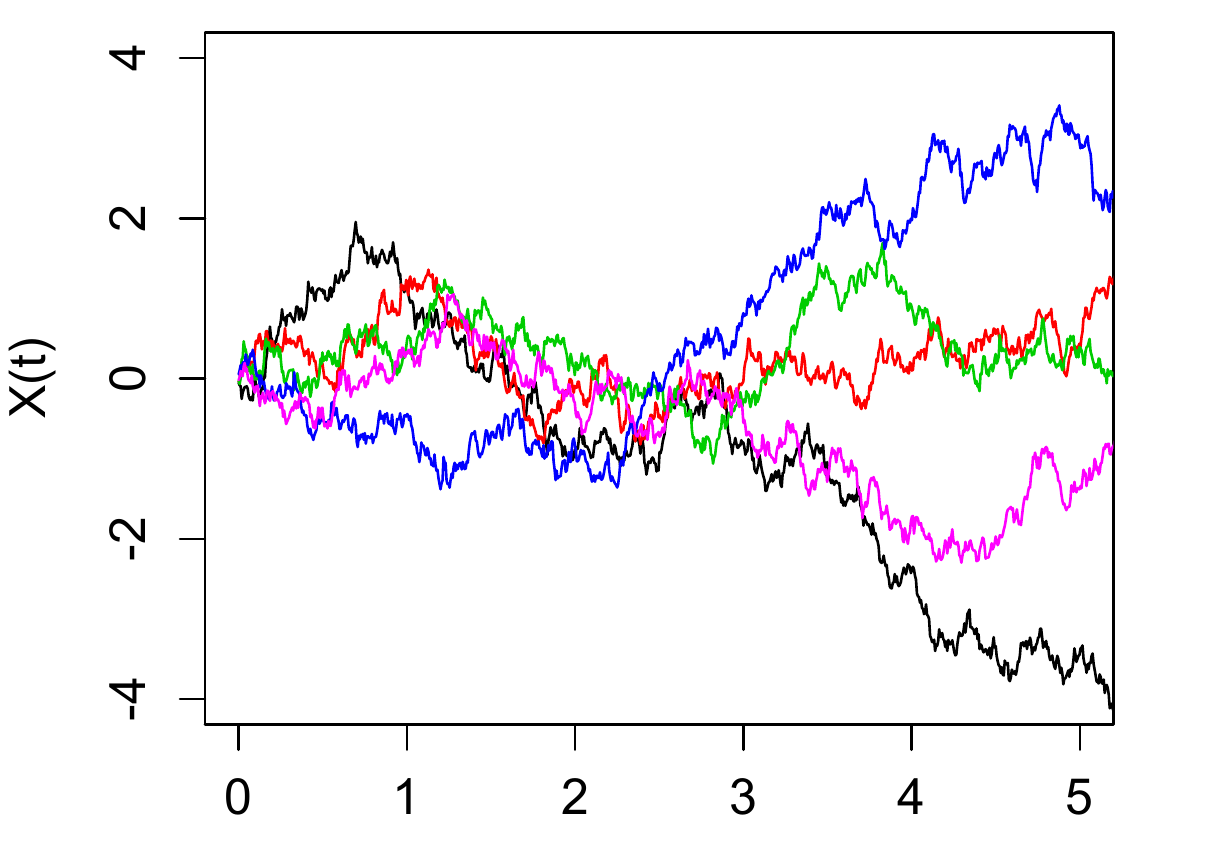}}
	\end{center}
\caption{Realizations of a Brownian motion with a diffusion coefficient of one half so that mean square displacement is equal to $t$.} \label{fig:sp4}
\end{figure}

\section{Fluctuations}
\label{sec:fluc}

We explained in \S\ref{markov-jump-proess} that  the fluctuations of a diffusion where the jump-rate satisfies \eqref{assump-jr} can take on a myriad of forms. Since the activity and variation of the sample path can be finite and infinite, four categories of particle fluctuations emerge. The activity refers to the number of jumps, or points of discontinuity, on a finite time interval. An infinite activity process contains a countable number of jumps and is associated with a jump-rate that is not integrable, i.e.,
\begin{align}
\int_{-\lambda}^\lambda \gamma(x,y)\, dy &= \infty \text{ for all } \lambda>0\,.\label{inf-act}
\end{align}
The variation of a one-dimensional function is the vertical distance traversed along the graph. The definition is somewhat more technical when the function is $n>1$ dimensional but embodies the same concept.

Figures~\ref{fig:sp1}-\ref{fig:sp4} display the four categories of particle fluctuations and provide a useful visual description.
Figures~\ref{fig:sp1} displays a finite activity, finite variation process induced by a compound Poisson process; the resemblance to a continuous-time Markov chain is apparent. Figure \ref{fig:sp4} displays Brownian motion, a ``zero activity'' process (since there are a finite number---zero---of jumps) of infinite variation.
The infinite activity processes displayed in Figures~\ref{fig:sp2}--\ref{fig:sp3} warrant further discussion since a countable number of jumps over a finite time interval cannot be computed but can be approximated.
\begin{subequations}\label{infinite-acivity}
Define the radial jump-rate kernel
\begin{align} \label{nu}
\nu(y) & := \frac{1}{m}\frac{1}{|y|^{1+\alpha}} \mathds{1}_{(0,\lambda)}(|y|)\,, \quad 0 < \alpha < 2\,,
\end{align}
and note that this kernel satisfies the condition \eqref{inf-act} and represents a truncated $\alpha$-stable L\'{e}vy process where the parameter $m$ controls the rate of diffusion and is used primarily so that the axes of Figures~\ref{fig:sp1}-\ref{fig:sp4} are aligned. In particular, when $0 < \alpha < 1$, then a finite variation process occurs because
\begin{align}
\sint{-\lambda}{\lambda}{|y| \nu(y)}{y} & < \infty\,,
\intertext{whereas when $1< \alpha < 2$, an infinite variation process occurs because $|y| \nu(y)$ is not integrable but}
\sint{-\lambda}{\lambda}{ y^2 \nu(y)}{y} & < \infty\,.
\end{align}

In order to simulate such a process, we employ the standard procedure to approximate $\nu$  with the integrable jump-rate
\begin{align} \label{nu-app}
		\nu_\varepsilon(y) & =
		\begin{cases}
			\ds \frac{1}{m}\frac{1}{|y|^{1+\alpha}} & \varepsilon < |y| < \lambda\,, \\
			\ds \frac{1}{m}\frac{1}{\varepsilon^{1+\alpha}} & |y| < \varepsilon\,,
		\end{cases}
\end{align}
with small $\varepsilon$.
Such an approximation enables us to realize a compound Poisson process that serves to approximate the infinite activity process.
\end{subequations}

\section{Nonlocal vector calculus} \label{sec:NLVC}

We briefly review the nonlocal vector calculus introduced in \cite{gule:10} and extensively developed in \cite{dglz:13}. The calculus will enable us to express the nonlocal diffusion equation in conservative form and demonstrate that the variational form of the equation is well-posed. 

Let $\balpha,\fb: \mathbb{R}^n \times \mathbb{R}^n \to \mathbb{R}^n$ where $\balpha$ is antisymmetric in $x$ and $y$, i.e., $ \balpha(x,y )=-\balpha(y,x )$. Define the nonlocal divergence
\begin{align} \label{def-D}
\mcD(\fb)(x) := \int_{\mathbb{R}^n} \big(\fb(x, y ) + \fb(y ,
x)\big) \cdot \balpha(x,y )\,dy\,.
\end{align}
The nonlocal divergence of a vector is a scalar and the definition implies that the antisymmetric part of $\fb$, i.e., $1/2\big(\fb(x,y )-\fb(y,x )\big) $ is annihilated.
If we suppose that $\fb$ is differentiable, then the choice
\begin{subequations} \label{dist-nd}
\begin{align}
\balpha(x,y) & = -\nabla_y \, \delta(y -x) \label{ddirac}
\intertext{and an integration by parts implies that}
\mcD(\fb)(x) & \equiv \nabla \cdot \fb(x,x)\,. \label{dist-div}
\end{align}
\end{subequations}
The antisymmetry of the integrand of \eqref{def-D} grants the nonlocal divergence theorem
\begin{subequations} \label{nl-div-thm}
\begin{align}
\int_{\Omega} \mcD(\fb)\,dx &= -\int_{\mathbb{R}^n\setminus\Omega} \mcD(\fb)\,dx \,.
\intertext{If we restrict $\balpha$ so that $\balpha(x,y)=\mathbf{0}$ for $x\in \Omega$ and $y \in \mathbb{R}^n\setminus\big(\Omega\cup \OmegaI\big)$ where the interaction domain $\OmegaI$ was defined in \S\ref{sec:process-cnfd-dm}, then}
\int_{\mathbb{R}^n\setminus\Omega} \mcD(\fb)\,dx & = \int_{\OmegaI} \mcD(\fb)\,dx\,. \label{nl-div}
\intertext{Both integrals represent the flux of the vector field into the region external to $\Omega$. The nonlocal divergence theorem is the nonlocal analogue of the classical divergence theorem}
\int_{\Omega} \nabla \cdot \fb \,dx & = \int_{\partial \Omega}  \fb \cdot \nb \,d^{n-1} x\,, \notag
\end{align}
\end{subequations}
where we abuse notation to suppose that $\fb$ is a vector field of only one variable. The relationship is clear when the relations \eqref{dist-nd} are invoked. The analogy between divergence theorems suggests that the orientation given by the unit normal $\nb$, or the sense of direction, is embodied by the antisymmetry of the kernel $\balpha$.

We now define the operator
\[
\big(\mcD^\ast u\big) (x,y ) := -\big(u(y )-u(x)\big)\,\balpha(x,y )
\]
so that given the density $u$, the function $\mcD^\ast u$ is a vector field of the points $x$ and $y$.  In a similar fashion to \eqref{dist-div}, the choice \eqref{ddirac} for the kernel $\balpha$ implies
\[
\int_{\mathbb{R}^n} \mcD^\ast u  \, dy = -\nabla u\,,
\]
so that we may conclude that the operator $\mcD^\ast$ can be seen as a nonlocal gradient. Because $\mcD^\ast u$ is an example of a vector field $\fb$, its nonlocal divergence can be determined. The identity
\begin{align}
\int_{\Omega\cup\OmegaI} v \, \mcD\,\fb\,dx = \int_{\Omega\cup\OmegaI}\int_{\Omega\cup\OmegaI} \mcD^\ast v \cdot \fb \,dy \, dx \,, \label{adj-ident}
\end{align}
is established by direct substitution and invoking the antisymmetry of $\balpha$. In words,  the identity states that  $\mcD^\ast$ is the adjoint operator for $\mcD$. This is analogous to the classical identity that the divergence $\nabla \cdot$ is the negative adjoint operator for the gradient $\nabla$.   Inserting $\fb = \mcD^\ast u$ into the previous identity and using linearity of the integral on the left-hand side grants a nonlocal Green's identity
\begin{align}
\int_{\Omega} v \mcD\,\mcD^\ast u\,dx &= \int_{\Omega\cup\OmegaI}\int_{\Omega\cup\OmegaI} \mcD^\ast v \cdot\mcD^\ast u \,dy \, dx - \int_{\OmegaI} v \, \mcD \, \mcD^\ast u\,dx\,,
\intertext{which is seen to be the nonlocal analogue of the conventional Green's identity}
\int_{\Omega} v \Delta u \,dx & = -\int_{\Omega} \nabla v \cdot \nabla u \,dx + \int_{\partial \Omega} v \big(\nabla u \cdot \nb\big) \, d^{n-1} x \,.
\end{align}
Let $v$ be a smooth, compactly supported function. Then, we can show that
\[
  \int_{\OmegaI} v \, \mcD\,\mcD^\ast u\,dx = \int_{\partial \Omega} v\, \big(\mcD^\ast u \cdot \mathbf{n} \big) \, d^{n-1} x\,;
\]
see \cite[\S5.1]{dglz:13} for details. This demonstrates in what sense the nonlocal identity is a generalization of the conventional identity by avoiding spatial derivatives and replacing surfaces by volumes for boundary and volume data, respectively. This more general Green's identity will allow us to demonstrate that the deterministic exit-time problem is associated with a general class of Markov jump processes is well-posed. 

\section{Nonlocal convection-diffusion equation} \label{sec:ncde}

We define the nonlocal convection-diffusion equation
  \[
  \left\{\begin{aligned}
  u_t &= \mcD \, \mathbf{f}\\
  \mathbf{f} & = \bmu \,u - \bthe \,\mcD^\ast u\,,
  \end{aligned}\right.
  \]
  where $\boldsymbol{\mu}: \mathbb{R}^n \times \mathbb{R}^n \to \mathbb{R}^n$ and $\bthe: \mathbb{R}^n \times \mathbb{R}^n \to \mathbb{R}^n\times \mathbb{R}^n$. The first and second lines represent a nonlocal Fick's first and second law of diffusion and are the analogues of the classical laws
  \[
  \left\{\begin{aligned}
  u_t &= -\nabla \cdot \mathbf{q}\\
  \mathbf{q} & = \mathbf{b} \,u + \nabla \mathbf{A} \, u\,.
  \end{aligned}\right.
  \]
Combining Fick's first and second laws leads to the nonlocal convection-diffusion equation
\begin{subequations} \label{nl-cdeq}
  \begin{align}
  u_t &= \mcD \, \big( \boldsymbol{\mu} \,u - \bthe \,\mcD^\ast u\big) \qquad \text{on } \Omega\cup\OmegaI\,. \label{nl-cdeq-a}
\intertext{This equation is simply a rewrite of the confined master equation \eqref{conf-ms-eq} with the identification}
 \gamma  & = \balpha\cdot\bthe\balpha - \bmu\cdot\balpha\,. \label{nl-cdeq-b}
 \end{align}
An immediate consequence for equation \eqref{nl-cdeq-a} is
 \begin{align}
\label{nl-cdeq-c}
  \begin{aligned}
   \frac{d}{dt} \int_{\Omega\cup\OmegaI} u \, dx & = \int_{\Omega\cup\OmegaI} \mcD \, \big( \boldsymbol{\mu} \,u - \bthe \,\mcD^\ast u\big) \, dx \\
& = -\int_{\OmegaI} \mcD \, \big( \boldsymbol{\mu} \,u - \bthe \,\mcD^\ast u\big) \,dx + \int_{\OmegaI} \mcD \, \big( \boldsymbol{\mu} \,u - \bthe \,\mcD^\ast u\big) \,dx \\
&= 0\,,
  \end{aligned} 
 \end{align}
\end{subequations}
where we invoked the nonlocal divergence theorem \eqref{nl-div-thm} for the second equality. In words, the
probability flux out of $\Omega$ into $\OmegaI$ is equal and opposite to the probability flux out of $\OmegaI$ into $\Omega$. This is an instance of the more general principle of action-reaction, where $\OmegaI$ and $\Omega$ can be replaced by regions that have no overlap and are separated by a finite distance. Hence the interaction is not restricted to regions that are in direct contact so leading to nonlocal diffusion.

The paper \cite[Theorem~2.1]{dglz:13} demonstrates that a well-formulated nonlocal balance law is given by the following four equivalent conditions:
\begin{enumerate}
\item
antisymmetry of $\psi(x,y):=u(y,t)\gamma(y,x) - u(x,t)\gamma(x,y)$,  the integrand of $\mcD \, \big( \boldsymbol{\mu} \,u - \bthe \,\mcD^\ast u\big)$;
\item  no self interaction, i.e., $\ds
\int_{\Omega} \mcD \, \big( \boldsymbol{\mu} \,u - \bthe \,\mcD^\ast u\big) \, dx =0$ for all $\Omega$;
\item action-reaction, i.e., $\ds \int_{\Omega} \int_{\Omega^\prime} \psi(x,y) \, dy \, dx + \int_{\Omega^\prime} \int_{\Omega} \psi(x,y) \, dy \, dx = 0$ for all $\Omega,\Omega^\prime$ that have no overlap;
    \item additivity, i.e., $\ds \int_{\Omega\cup\Omega^\prime } \mcD \, \big( \boldsymbol{\mu} \,u - \bthe \,\mcD^\ast u\big) \, dx = \int_{\Omega} \mcD \, \big( \boldsymbol{\mu} \,u - \bthe \,\mcD^\ast u\big) \, dx + \int_{\Omega^\prime } \mcD \, \big( \boldsymbol{\mu} \,u - \bthe \,\mcD^\ast u\big) \, dx$ for all $\Omega,\Omega^\prime$ that have no overlap.
\end{enumerate}

\section{Backward and forward Kolmogorov equations}\label{sec:kol}

A consequence of the identity \eqref{nl-cdeq-c} implies that the equation \eqref{nl-cdeq-a} is formally understood as a forward Kolmogorov, or Fokker-Planck, equation because the transition measure for the jump diffusion is evolved. Convention denotes the corresponding operator via its action on the density
\begin{subequations}
  \begin{align}
\mathcal{A}^\ast u & :=  \mcD \, \big( \boldsymbol{\mu} \,u\big) -  \mcD \,  \bthe \,\mcD^\ast u \,, \label{op}
\intertext{where $\mathcal{A}^\ast$ denotes the adjoint operator for $\mathcal{A}$. In order to consider the mean exit-time problem, the nonlocal backward Kolmogorov equation }
\label{nl-bkeq}  u_t &= \mathcal{A} \, u 
\intertext{is also helpful. 
An explicit expression for $\mathcal{A}$, the generator for the Markov process, is possible by exploiting the nonlocal vector calculus. To simplify matters, we assume, without loss of generality that
$
\bthe(x,y)=\bthe^T(x,y)= \bthe(y,x)\,,
$
so that by the adjoint identity \eqref{adj-ident} with $\fb(x,y) = \bmu(x,y)\,u(x)$, the requisite operator is}
\big(\mathcal{A} \, u\big) (x) &:= \sint{\Omega\cup\OmegaI}{}{\big(\bmu \cdot \mcD^\ast u\big)(x,y)}{y} - \big(\mcD \, \bthe \,\mcD^\ast\big) u(x)
\,. \label{adj-op}
\end{align}
\end{subequations}
When the jump-rate satisfies $\gamma(x,y)=\gamma(y,x)$, or equivalently by \eqref{nl-cdeq-b}, when $\bmu(x,y)=-\bmu(y,x)$, then $\mathcal{A}=\mathcal{A}^\ast$. 

\section{Exit-time problem}

Let $X_t$ be a finite-range Markov jump process for the confined master equation \eqref{conf-ms-eq} conditioned on $X_0=x\in\Omega$ and let a generalized exit-time random variable be given by
\begin{align}
T_x & := \inf\{ t : X_t \in \Omega_d \subseteq \OmegaI \}\,,
\end{align}
where $\emptyset \subseteq \Omega_d \subseteq \OmegaI$. We denote the exit-time generalized since the particle cannot exit to the region of $\OmegaI$ not in $\Omega_d$ (except when $\Omega_d \equiv \OmegaI$).

In contrast to the classical exit-time problem for Brownian motion, we cannot expect the finite-range jump process $X_t$ to hit the boundary $\partial \Omega$ upon departing $\Omega$. In fact, with probability almost surely one, the process never encounters the boundary. Instead, $X_t$ departs to a location in the interaction domain $\OmegaI$.
The density of particles that have not yet exited $\Omega$ to $\Omega_d $ evolves according to the nonlocal convection-diffusion, or Fokker-Planck, equation:
\begin{equation} \label{nonlocal-exit}
	\left\{\begin{aligned}
			\ds u_t(x,t) & = \big(\mathcal{A}^\ast u \, \mathds{1}_{\Omega\cup \Omega_d}\big) (x,t) \\
			& = \sint{\Omega \cup \Omega_d}{}{\big( u(y,t)\gamma(y,x) - u(x,t)\gamma(x,y) \big)}{y},  & x &\in {\Omega} \\
			\ds u(x,t) & = 0, & x &\in \Omega_d \\
			\ds u(x,0) & = u_0(x),  & x & \in \Omega\,,
		\end{aligned} \right.
	\end{equation}
where we used \eqref{indicator}, \eqref{nl-cdeq-b} and \eqref{op} for the second equality, and without loss of generality we suppose that $u_0$ is a probability density over $\Omega$.
We define the constraint on the density over $\Omega_d$ to be a volume-constraint, the generalization of a boundary condition, or committing a semantic abuse of terminology, a nonlocal boundary condition. The volume constraint is necessary because the particle may jump out of $\Omega$ into $\OmegaI$. As we will review in \S\ref{sec:anal-nlde}, the volume constraint is crucial in showing that the nonlocal convection-diffusion equations is well-posed. Moreover, the volume-constraint enables us to link the finite-range jump process with a deterministic equation and corresponds to how Monte Carlo realizations are implemented.

\begin{figure}[htbp]
	\begin{center}
	\resizebox{0.7\columnwidth}{!}{\includegraphics{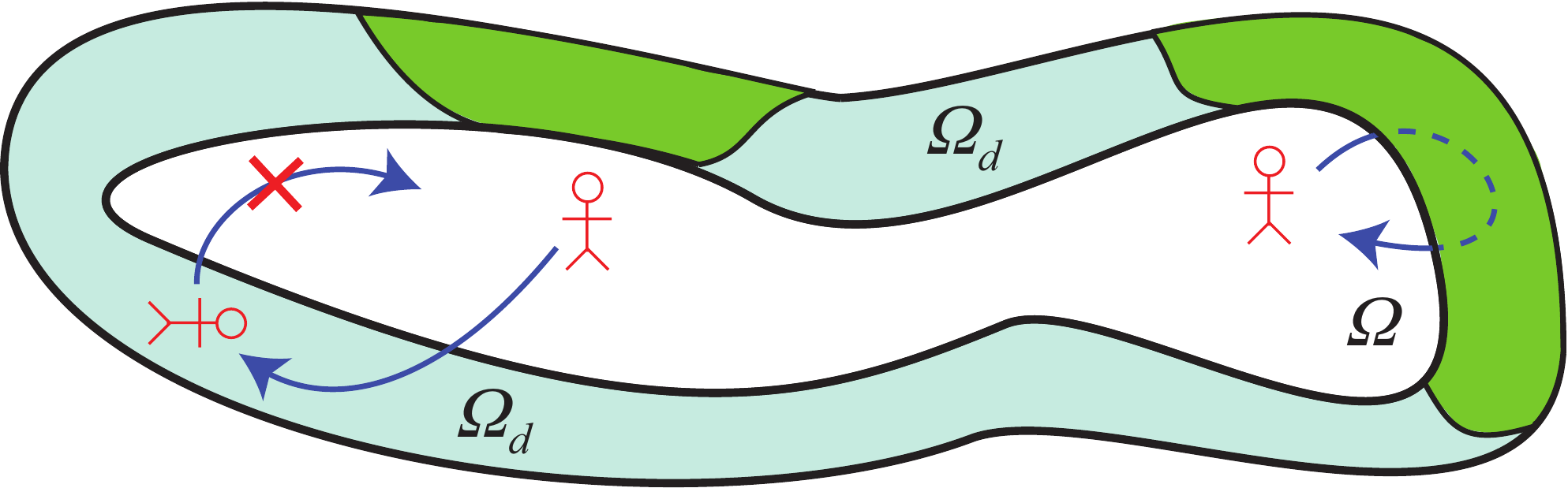}}
	\end{center}
\caption{Absorbed/censored jump process depicted by the possible steps taken by the random walker stickman.} \label{fig:mixed}
\end{figure}

Two cases are of special interest.
\begin{description}
  \item[$\Omega_d \equiv \emptyset$:]
   The system \eqref{nonlocal-exit} is then a nonlocal analogue of the classical Fokker-Planck with a homogenous Neumann boundary condition. The corresponding jump process models a particle confined to the region $\Omega$ and is an instance of a censored process; see the paper \cite{bobc:03} for a mathematical discussion.
  \item[$\Omega_d \equiv \OmegaI$:]
     The system \eqref{nonlocal-exit} is then a nonlocal analogue of the classical Fokker-Planck with a homogenous Dirichlet boundary condition.
     The corresponding jump process models the exit-time of a particle conditioned on $X_0\in\Omega$.
\end{description}
The general case corresponds to a system for a mixed absorbed/censored process and is the analogue of the classical Fokker-Planck with a mixed Dirichlet/Neumann boundary condition.  Figure~\ref{fig:mixed} depicts this general case and this section provides the infrastructure needed to extend the analyses of the narrow escape problem~\cite{hosc:14} beyond classical diffusion to Markov jump processes.

The homogenous Dirichlet volume-constraint is explicit in the system \eqref{nonlocal-exit}; where is the homogenous Neumann volume-constraint? Suppose that instead of a homogenous condition, we consider a given source density $g$. Then we may augment \eqref{nonlocal-exit} with the condition
\[
\sint{\OmegaI\setminus\Omega_d}{}{\big( u(y,t)\gamma(y,x) - u(x,t)\gamma(x,y) \big)}{y} = g(x) \qquad x \in \Omega\,,
\]
to obtain a modified system. If we specify $g \equiv 0$ over $\Omega$ then we see that the homogenous Neumann volume-constraint is also explicitly given by the formulation of the system \eqref{nonlocal-exit}.

How does the deterministic system \eqref{nonlocal-exit} conserve the probability of the location of the particle? The derivation provides the insight needed to properly understand the distinction between the nonlocal and classical systems evolving the density. 
By invoking the volume-constraint, we have
\begin{align*}
\frac{d}{dt} \int_{\Omega} u(x,t) \, dx & = \int_{\Omega} \sint{\Omega \cup \Omega_d}{}{\big( u(y,t)\gamma(y,x) - u(x,t)\gamma(x,y) \big)}{y}\, dx \\
& = \int_{\Omega} \sint{\Omega_d}{}{\big( u(y,t)\gamma(y,x) - u(x,t)\gamma(x,y) \big)}{y}\, dx \\
& = -\int_{\Omega} \sint{\Omega_d}{}{u(x,t)\gamma(x,y)}{y}\, dx < 0\,,
\end{align*}
where the second and third equalities follow from the antisymmetry of the integrand and applying the volume constraint a second time, respectively. By the nonlocal divergence theorem, we may also conclude that the last iterated integral is the nonlocal flux of probability out of $\Omega$ into $\Omega_d$. Therefore the probability is conserved
over $\Omega \cup \Omega_d$, i.e., the particle is located either in $\Omega$ or has exited to $\Omega_d$.

We now provide a probabilistic interpretation. Because the probability that the particle remains in $\Omega$ is given by
\begin{align*}
\text{Prob}(T_x>t) &= \int_\Omega u(y,t)\, dy\,,
\intertext{then the exit-time distribution is given by}
\text{Prob}(T_x\leq t) & = 1 - \int_\Omega u(y,t)\, dy \\
& = \int_0^t \int_{\Omega} \sint{\Omega_d}{}{u(y,s)\gamma(y,z)}{y}\, dz \, ds\,,
\end{align*}
where the last equality is the time integrated probability flux.
In words, the rate of change of the probability of the particle exiting to $\Omega_d$ increases to one so that as time increases the density of particles located in $\Omega$ decreases to zero. We also see the effect of decreasing the size of $\Omega_d$ is to delay the time required for the probability of a particle to exit to $\Omega_d$. And when $\Omega_d \equiv \emptyset$, then the probability of a particle to exiting $\Omega$ is zero for all time.

The mean exit-time
\[
\langle T_x \rangle : = \int_0^\infty \int_\Omega u(x,t) \, dx \, dt\,,
\]
then, in direct analogy with the classical exit time problem, see, e.g., the paper \cite{nkms:90}, is given by the solution of the steady-state volume-constrained problem
\begin{equation} \label{nonlocal-mean-exit}
	\left\{\begin{aligned}
			\mathcal{A} \, \langle T_x \rangle & = -1  & x &\in {\Omega} \\
			 \langle T_x \rangle & = 0, & x &\in \Omega_d\,\, (\neq \emptyset)\,,
		\end{aligned} \right.
\end{equation}
where the operator $\mathcal{A} $ is given by \eqref{adj-op}. An equation for the remaining moments of the exit-time random variable $T_x$ may also be determined via repeated integration, in an analogous fashion to the classical case; see, e.g., \cite{gard:04}.

An important question is whether the exit-time distribution and moments are finite. Given sufficient conditions on the kernel $\gamma$, the distribution and moments are indeed finite; see \S\ref{sec:anal-nlde} for a discussion.

The above analysis allows us to consider far more interesting problems than that depicted by Figure~\ref{fig:mixed}; for instance, consider the Figure~\ref{fig:vc2}. Such an exit-time analysis is simply not possible with Brownian motion since $\Omega$ is the union of a disconnected set of regions---a particle undergoing Brownian motion cannot jump outside of $\Omega$. The report \cite{bule:13} provides an example of such an analysis.

\section{Analysis and approximation of nonlocal diffusion equation}
\label{sec:anal-nlde}

The  goal of this section is to briefly review analyses demonstrating that the volume-constrained nonlocal diffusion equation is well-posed. For simplicity, we suppose that the interaction domain $ \OmegaI \equiv \Omega_d$. 
A crucial aspect of the analysis is that a nonlocal variational formulation based upon the nonlocal vector calculus introduced in \S\ref{sec:NLVC} is exploited.

The analysis proceeds by first demonstrating that the steady-state problem
\begin{align}
&\left\{\begin{aligned}
\mathcal{A}^\ast u & = b & \text{on } & \Omega \\
u & = 0 & \text{on } & \OmegaI  \,,
\end{aligned}\right. \label{stysteq}
\end{align}
where $\mathcal{A}^\ast$ was defined in \eqref{op}, is well-posed. Then, standard results are invoked to demonstrate that the time dependent nonlocal diffusion equation is well-posed. 
This sets the stage for the development of stable numerical methods for the discretization of \eqref{nonlocal-exit} that offer an alternative to averaging the results of realizations of the jump process; see \cite{bule:11,bule:13,duhl:14} for numerical results demonstrating the consistency of the numerical solutions of the volume constrained diffusion equation and Monte Carlo simulations.

We now review the analysis of the steady-state problem. Given a Hilbert space $V$ with the inner product between elements $u,v$
\begin{align}
(u,v)_V & :=\int_{\Omega\cup\OmegaI}\int_{\Omega\cup\OmegaI} \big( u\,v + \mcD^\ast u \cdot \mcD^\ast v \big) \,dy \, dx\,,\notag
\intertext{the variational problem is: Find $u \in V$ such that}
a(u,v) & = \int_{\Omega} u\, b \, dx \qquad \forall \, v \in V \label{auv=b}
\intertext{where the bilinear form $a(\cdot,\cdot)$ is given by}
  a(u,v) & = \int_{\Omega\cup\OmegaI}\int_{\Omega\cup\OmegaI} \mcD^\ast u \cdot \bthe \mcD^\ast v \,dy \, dx - \int_{\Omega} \mcD\big(\boldsymbol{\mu} u \big) \, v\, dx \notag \\
& = \int_{\Omega\cup\OmegaI}\int_{\Omega\cup\OmegaI} \big(u(y)-u(x)\big) \big(v(y)-v(x)\big) \, \big(\balpha\cdot \bthe \balpha\big)(x,y)\,dy \, dx \notag \\
& \, \,\,\,- \int_{\Omega}\int_{\Omega\cup\OmegaI} \big( \bmu(x,y)u(x) + \bmu(y,x)u(y)\big)\,v(x) \,dy \, dx\,, \notag
\end{align}
and an expression for $\gamma$ in terms of $\balpha, \bthe, \bmu$ was given in \eqref{nl-cdeq}.
The Lax-Milgram theorem then provides sufficient conditions that when satisfied, demonstrate that \eqref{auv=b} has a unique density $u$, i.e., the system \eqref{stysteq} has a weak solution. 
The Hilbert space $V$ is identified with a volume-constrained subspace of square integrable functions or a fractional Sobolev space given conditions on the integrability of the jump-rate $\gamma$. The latter space contains the densities corresponding to infinite activity processes.
Continuous dependence upon the data implies the energy is bounded by the data, i.e.,
\begin{align}
(u,u)_V  \leq C \int_\Omega b^2 \, dx\,; \notag
\end{align}
see \cite{duhl:14,ddgl:14} for details and further discussion.
If $\gamma(x,y)=\gamma(y,x)$ then the bilinear form $a(\cdot,\cdot)$ is symmetric and the variational problem \eqref{auv=b} is the Euler-Lagrange equation for the minimization problem
\[
    \min_{V}  \frac{1}{2} a(u,v) - \int_{\Omega} u\, b \, dx\,,
\]
this is the case considered in \cite{dglz:12}. To the best of our knowledge, this latter paper was the first to demonstrate that the deterministic exit-time problem for an infinite activity and finite variation jump diffusion is well-posed. The paper \cite{degu:13} exploits the analysis in \cite{dglz:12} to demonstrate how the truncated fractional Laplacian converges to the fractional Laplacian as the length $\lambda$ of the finite-range increases.

One of the hypothesis associated with the Lax-Milgram theorem is to show that the bilinear form $a(\cdot,\cdot)$ is coercive. Equivalently, coercivity means that
\begin{align*}
	\sigma = \inf_{v \in V} \frac{a(v,v)}{(u,v)_V} \,. \label{dirichlet_eigenvalue}
\end{align*}
is a positive number. Coercivity can then be used to demonstrate that the moments of the exit-time are finite; see \cite{bule:13,ddgl:14} for a demonstration. 

\section*{Acknowledgements}

We thank Professors Scott McKinley (University of Florida) and Renming Song (University Illinois) for helpful discussions, and Professor Michael Mascagni of Florida State University for a careful reading. 

\bibliographystyle{abbrv}
\bibliography{nonlocal-diffusion}

%
%

\end{document}